# On the Developable Mannheim Offsets of Spacelike Ruled Surfaces

**Mehmet Önder[1], H. Hüseyin Uğurlu[2]**
[1]*Celal Bayar University, Faculty of Arts and Sciences, Department of Mathematics, Muradiye Campus, 45047,Muradiye, Manisa, Turkey.*
[2]*Gazi University, Gazi Faculty of Education, Department of Secondary Education Science and Mathematics Teaching,Mathematics Teaching Program, Ankara, Turkey*
E-mails: mehmet.onder@cbu.edu.tr, hugurlu@gazi.edu.tr

**Abstract**
In this paper, using the classifications of timelike and spacelike ruled surfaces, we define and study the Mannheim offsets of spacelike ruled surfaces in Minkowski 3-space. We give the conditions for spacelike offset surfaces to be developable.

## 1. Introduction

Ruled surfaces are the surfaces which are generated by moving a straight line continuously in the space and are one of the most important topics of differential geometry. A ruled surface can always be described (at least locally) as the set of points swept by a moving straight line. These kinds of surfaces are used in many areas of sciences such as Computer Aided Geometric Design (CAGD), mathematical physics, moving geometry and kinematics for modeling the problems and engineering. So, many geometers and engineers have studied on ruled surfaces in different spaces such as Euclidean space and Minkowski space [1,2,3,5,6,11,12,14]. They have investigated and obtained many properties of the ruled surfaces in these spaces.

After the discovery of preconstraint concrete in 1930, ruled surfaces have had an important role in architectural construction. These surfaces have been used to construct the water-towers, chimney-pieces, roofs, and spiral staircases. Eero Saarinen (1910-1961) used ruled surfaces in his buildings at Yale and M.I.T. Furthermore, builder Félix Candela has made extensive use of cylinders and the most familiar ruled surfaces [4].

Differently from Euclidean space, in the Minkowski space a ruled surface has different properties according to the Lorentzian casual characters of the lines and curves of the surface. The classifications of the ruled surfaces in the Minkowski 3-space $IR_1^3$ are given by Kim and Yoon [5]. Using the classification of the ruled surfaces with non-null frame vectors and their non-null derivatives, Önder and Uğurlu have given the Frenet frames, invariants and instantaneous rotation vectors of the Frenet frames of the timelike and spacelike ruled surfaces in Minkowski 3-space $IR_1^3$ [10,15].

Furthermore, in the plane, while a curve $\alpha$ rolls on a straight line, the center of curvature of its point of contact describes a curve $\beta$ which is called Mannheim of $\alpha$. Mannheim partner curves in three dimensional Euclidean 3-space and Minkowski 3-space $IR_1^3$ are studied by Liu and Wang [7,17]. They have given the definition of Mannheim partner curves as follows: Let $C$ and $C^*$ be two space curves. $C$ is said to be a Mannheim partner curve of $C^*$ if there exists a one to one correspondence between their points such that the binormal vector of $C$ is the principal normal vector of $C^*$. They showed that $C$ is Mannheim partner curve of $C^*$ if and only if the following equality holds

$$\frac{d\tau}{ds} = \frac{\kappa}{\lambda}(1+\lambda^2\tau^2),$$

where $\kappa$ and $\tau$ are the curvature and the torsion of the curve $C$, respectively, and $\lambda$ is a nonzero constant [17].

Ravani and Ku have given the correspondence of the notion of Bertrand curves to the ruled surfaces and called Bertrand offsets of ruled surfaces [18]. Similarly, by considering the notion of Mannheim partner curves, Orbay, Kasap and Aydemir have defined the Mannheim offsets of the ruled surfaces in 3-dimensional Euclidean space $E^3$ [9].

In this paper, by considering the classifications of the ruled surfaces in Minkowski 3-space $IR_1^3$, we give the Mannheim offsets of spacelike ruled surfaces and we obtain the conditions for spacelike offset surfaces to be developable.

## 2. Preliminaries

The Minkowski 3-space $IR_1^3$ is the real vector space $IR^3$ provided with the standard flat metric given by
$$\langle , \rangle = -dx_1^2 + dx_2^2 + dx_3^2$$
where $(x_1, x_2, x_3)$ is a standard rectangular coordinate system of $IR_1^3$. An arbitrary vector $\vec{v} = (v_1, v_2, v_3)$ in $IR_1^3$ can have one of three Lorentzian causal characters; it can be spacelike if $\langle \vec{v}, \vec{v} \rangle > 0$ or $\vec{v} = 0$, timelike if $\langle \vec{v}, \vec{v} \rangle < 0$ and null (lightlike) if $\langle \vec{v}, \vec{v} \rangle = 0$ and $\vec{v} \neq 0$ [8]. Similarly, an arbitrary curve $\vec{\alpha} = \vec{\alpha}(s)$ is spacelike, timelike or null (lightlike), if all of its velocity vectors $\vec{\alpha}'(s)$ are spacelike, timelike or null (lightlike), respectively. We say that a timelike vector is future pointing or past pointing if the first compound of the vector is positive or negative, respectively. Let $\vec{v} = (v_1, v_2, v_3) \in IR_1^3$ be a curve. Then $\|\vec{v}\| = \sqrt{|\langle \vec{v}, \vec{v} \rangle|}$ is called the norm of the vector $\vec{v}$.

For any vectors $\vec{x} = (x_1, x_2, x_3)$ and $\vec{y} = (y_1, y_2, y_3)$ in $IR_1^3$, in the meaning of Lorentz vector product of $\vec{x}$ and $\vec{y}$ is defined by
$$\vec{x} \times \vec{y} = \begin{vmatrix} e_1 & -e_2 & -e_3 \\ x_1 & x_2 & x_3 \\ y_1 & y_2 & y_3 \end{vmatrix} = \left( \begin{vmatrix} x_2 & x_3 \\ y_2 & y_3 \end{vmatrix}, \begin{vmatrix} x_1 & x_3 \\ y_1 & y_3 \end{vmatrix}, -\begin{vmatrix} x_1 & x_2 \\ y_1 & y_2 \end{vmatrix} \right).$$

The Lorentzian sphere and hyperbolic sphere of radius $r$ and center 0 in $IR_1^3$ are given by
$$S_1^2 = \{\vec{x} = (x_1, x_2, x_3) \in IR_1^3 : \langle \vec{x}, \vec{x} \rangle = r^2\},$$
and
$$H_0^2 = \{\vec{x} = (x_1, x_2, x_3) \in IR_1^3 : \langle \vec{x}, \vec{x} \rangle = -r^2\},$$
respectively [16].

**Definition 2.1.** ([13]) *i) Central angle:* Let $\vec{x}$ and $\vec{y}$ be spacelike vectors in $IR_1^3$ that span a timelike vector subspace. Then there is a unique real number $\theta \geq 0$ such that $|<\vec{x}, \vec{y}>| = \|\vec{x}\|\|\vec{y}\|\cosh\theta$. This number is called the *central angle* between the vectors $\vec{x}$ and $\vec{y}$.

*ii) Lorentzian timelike angle:* Let $\vec{x}$ be a spacelike vector and $\vec{y}$ be a timelike vector in $IR_1^3$. Then there is a unique real number $\theta \geq 0$ such that $|<\vec{x}, \vec{y}>| = \|\vec{x}\|\|\vec{y}\|\sinh\theta$. This number is called the *Lorentzian timelike angle* between the vectors $\vec{x}$ and $\vec{y}$.

**Definition 2.2.** A surface in the Minkowski 3-space $IR_1^3$ is called a timelike surface if the induced metric on the surface is a Lorentz metric and is called a spacelike surface if the induced metric on the surface is a positive definite Riemannian metric, i.e., the normal vector on the spacelike (timelike) surface is a timelike (spacelike) vector [16].

## 3. Ruled Surfaces in the Minkowski 3-space

Let $I$ be an open interval in the real line $IR$. Let $\vec{k} = \vec{k}(s)$ be a curve in $IR_1^3$ defined on $I$ and $\vec{q} = \vec{q}(s)$ be a unit direction vector of an oriented line in $IR_1^3$. Then we have the following parametrization for a ruled surface $N$

$$\vec{\varphi}(s, v) = \vec{k}(s) + v\vec{q}(s). \tag{1}$$

The parametric $s$-curve of this surface is a straight line of the surface which is called ruling. For $v = 0$, the parametric $v$-curve of this surface is $\vec{k} = \vec{k}(s)$ which is called base curve or generating curve of the surface. In particular, if the direction of $\vec{q}$ is constant, the ruled surface is said to be cylindrical, and non-cylindrical otherwise.

The striction point on a ruled surface $N$ is the foot of the common normal between two consecutive rulings. The set of the striction points constitute a curve $\vec{c} = \vec{c}(s)$ lying on the ruled surface and is called striction curve. The parametrization of the striction curve $\vec{c} = \vec{c}(s)$ on a ruled surface is given by

$$\vec{c}(s) = \vec{k}(s) - \frac{\langle d\vec{q}, d\vec{k} \rangle}{\langle d\vec{q}, d\vec{q} \rangle} \vec{q}. \tag{2}$$

[10,15]. So that, the base curve of the ruled surface is its striction curve if and only if $\langle d\vec{q}, d\vec{k} \rangle = 0$. Furthermore, the generator $\vec{q}$ of a developable ruled surface is tangent of its striction curve.

The distribution parameter (or drall) of the ruled surface in (1) is given as

$$\delta_\varphi = \frac{\left| d\vec{k}, \vec{q}, d\vec{q} \right|}{\langle d\vec{q}, d\vec{q} \rangle} \tag{3}$$

If $\left| d\vec{k}, \vec{q}, d\vec{q} \right| = 0$, then the normal vectors are collinear at all points of the same ruling and at the nonsingular points of the surface $N$, the tangent planes are identical. We then say that the tangent plane contacts the surface along a ruling. Such a ruling is called a *torsal* ruling. If $\left| d\vec{k}, \vec{q}, d\vec{q} \right| \neq 0$, then the tangent planes of the surface $N$ are distinct at all points of the same ruling which is called nontorsal. In Minkowski 3-space a ruled surface whose all rulings are torsal is called a developable ruled surface. The remaining ruled surfaces are called skew ruled surfaces. Then it is clear that in Minkowski 3-space, a ruled surface is developable if and only if at all its points the distribution parameter $\delta_\varphi = 0$ (For details see [10,15]).

For the unit normal vector $\vec{m}$ of the timelike ruled surface $N$ we have $\vec{m} = \frac{\vec{\varphi}_s \times \vec{\varphi}_v}{\|\vec{\varphi}_s \times \vec{\varphi}_v\|}$.

So, at the points of a nontorsal ruling $s = s_1$ we have

$$\vec{a} = \lim_{v \to \infty} \vec{m}(s_1, v) = \frac{d\vec{q} \times \vec{q}}{\|d\vec{q}\|}.$$

The plane of the timelike ruled surface $N$ which passes through its ruling $s_1$ and is perpendicular to the vector $\vec{a}$ is called the *asymptotic plane* $\alpha$. The tangent plane $\gamma$ passing through the ruling $s_1$ which is perpendicular to the asymptotic plane $\alpha$ is called the *central plane*. Its point of contact $C$ is *central point* of the ruling $s_1$. The straight lines which pass through point $C$ and are perpendicular to the planes $\alpha$ and $\gamma$ are called the *central tangent* and *central normal*, respectively.

Using the perpendicularly of the vectors $\vec{q}, d\vec{q}$ and the vector $\vec{a}$, representation of the unit vector $\vec{h}$ of the central normal is given by $\vec{h} = \frac{d\vec{q}}{\|d\vec{q}\|}$.

The orthonormal system $\{C; \vec{q}, \vec{h}, \vec{a}\}$ is called Frenet frame of the ruled surfaces where $C$ is the striction point [10,15].

Let now consider the ruled surface $N$. According to the Lorentzian characters of ruling and central normal, we can give the following classifications of the ruled surface $N$ as follows;

**i)** If the central normal vector $\vec{h}$ is spacelike and $\vec{q}$ is timelike, then the ruled surface $N$ is said to be of type $N_-$.

**ii)** If the central normal vector $\vec{h}$ and the ruling $\vec{q}$ are both spacelike, then the ruled surface $N$ is said to be of type $N_+$.

**iii)** If the central normal vector $\vec{h}$ is timelike, the ruling $\vec{q}$ and its derivative $d\vec{q}/ds$ are spacelike, then the ruled surface $N$ is said to be of type $N_\times$ [10,15].

The ruled surfaces of type $N_+$ and $N_-$ are clearly timelike and the ruled surface of type $N_\times$ is spacelike. By using these classifications the parametrization of the ruled surface $N$ can be given as follows,

$$\varphi(s,v) = \vec{k}(s) + v\vec{q}(s), \tag{4}$$

where $\langle \vec{h}, \vec{h} \rangle = \varepsilon_1 (=\pm 1)$, $\langle \vec{q}, \vec{q} \rangle = \varepsilon_2 (=\pm 1)$.

The set of all bound vectors $\vec{q}(s)$ at the point O constitutes the *directing cone* of the ruled surface $N$. If $\varepsilon_2 = -1$ (resp. $\varepsilon_2 = 1$), the end points of the vectors $\vec{q}(s)$ drive a spherical spacelike (resp. spacelike or timelike) curve $k_1$ on hyperbolic unit sphere $H_0^2$ (resp. on Lorentzian unit sphere $S_1^2$), called the *hyperbolic (resp. Lorentzian) spherical image* of the ruled surface $N$, whose arc is denoted by $s_1$.

For the Frenet vectors $\vec{q}, \vec{h}$ and $\vec{a}$ we have the following Frenet frames of ruled surface $M$:

**i)** If the ruled surface $N$ is timelike ruled surfaces of the type $N_+$ and $N_-$ then we have

$$\begin{bmatrix} d\vec{q}/ds_1 \\ d\vec{h}/ds_1 \\ d\vec{a}/ds_1 \end{bmatrix} = \begin{bmatrix} 0 & 1 & 0 \\ -\varepsilon_2 & 0 & \kappa \\ 0 & \varepsilon_2\kappa & 0 \end{bmatrix} \begin{bmatrix} \vec{q} \\ \vec{h} \\ \vec{a} \end{bmatrix}. \tag{5}$$

Darboux vector of the Frenet frame $\{O; \vec{q}, \vec{h}, \vec{a}\}$ can be given by $\vec{w}_1 = \varepsilon_2\kappa\vec{q} - \vec{a}$. Thus, for the derivatives in (5) we can write

$$d\vec{q}/ds_1 = \vec{w}_1 \times \vec{q}, \quad d\vec{h}/ds_1 = \vec{w}_1 \times \vec{h}, \quad d\vec{a}/ds_1 = \vec{w}_1 \times \vec{a},$$

and also we have

$$\vec{q} \times \vec{h} = \varepsilon_2\vec{a}, \quad \vec{h} \times \vec{a} = -\varepsilon_2\vec{q}, \quad \vec{a} \times \vec{q} = -\vec{h}. \tag{6}$$

[See 10].

**ii)** If the ruled surface $N$ is spacelike ruled surface of the type $N_\times$ then we have

$$\begin{bmatrix} d\vec{q}/ds_1 \\ d\vec{h}/ds_1 \\ d\vec{a}/ds_1 \end{bmatrix} = \begin{bmatrix} 0 & 1 & 0 \\ 1 & 0 & \kappa \\ 0 & \kappa & 0 \end{bmatrix} \begin{bmatrix} \vec{q} \\ \vec{h} \\ \vec{a} \end{bmatrix}. \tag{7}$$

Darboux vector of this frame is $\vec{w}_1 = -\kappa\vec{q} + \vec{a}$. Then the derivatives of the vectors of Frenet frame in (7) can be given by

$$d\vec{q}/ds_1 = \vec{w}_1 \times \vec{q}, \quad d\vec{h}/ds_1 = \vec{w}_1 \times \vec{h}, \quad d\vec{a}/ds_1 = \vec{w}_1 \times \vec{a}$$

and also we have

$$\vec{q} \times \vec{h} = -\vec{a}, \quad \vec{h} \times \vec{a} = -\vec{q}, \quad \vec{a} \times \vec{q} = \vec{h}. \tag{8}$$

[See 15].

In the equations (5) and (7), $s_1$ is the arc of generating curve $k_1$ and $\kappa = \dfrac{ds_3}{ds_1} = \left\| \dfrac{d\vec{a}}{ds} \right\|$ is conical curvature of the directing cone where $s_3$ is the arc of the spherical curve $k_3$ circumscribed by the bound vector $\vec{a}$ at the point $O$ [10,15].

## 4. Mannheim Offsets of the Spacelike Ruled Surfaces in Minkowski 3-space

Assume that $\varphi$ and $\varphi^*$ be two ruled surfaces in the Minkowski 3-space $IR_1^3$ with the parametrizations

$$\begin{cases} \varphi(s,v) = \vec{c}(s) + v\vec{q}(s), \\ \varphi^*(s,v) = \vec{c}^*(s) + v\vec{q}^*(s), \end{cases} \tag{9}$$

respectively, where $(\vec{c})$ (resp. $(\vec{c}^*)$) is striction curve of the ruled surfaces $\varphi$ (resp. $\varphi^*$). Let Frenet frames with non-null vectors of the ruled surfaces $\varphi$ and $\varphi^*$ be $\{\vec{q}, \vec{h}, \vec{a}\}$ and $\{\vec{q}^*, \vec{h}^*, \vec{a}^*\}$, respectively. Let $\varphi$ be a spacelike ruled surface. The ruled surface $\varphi^*$ is said to be Mannheim offset of the spacelike ruled surface $\varphi$ if there exists a one to one correspondence between their rulings such that the asymptotic normal of $\varphi$ is the central normal of $\varphi^*$. In this case, $(\varphi, \varphi^*)$ is called a pair of Mannheim ruled surface. By definition we can write

$$\vec{h}^* = \vec{a} \tag{10}$$

and so that by considering the classifications of the ruled surfaces it is easily seen that the Mannheim offsets of a spacelike ruled surface is a timelike ruled surfaces of the type $N_-$ or $N_+$. Then according to the Definition 2.1 and by using (10) we can give the following cases:

**Case 1.** If the Mannheim offset $\varphi^*$ of $\varphi$ is a timelike ruled surface of the type $N_-$ then we have

$$\begin{pmatrix} \vec{q}^* \\ \vec{h}^* \\ \vec{a}^* \end{pmatrix} = \begin{pmatrix} \sinh\theta & \cosh\theta & 0 \\ 0 & 0 & 1 \\ \cosh\theta & \sinh\theta & 0 \end{pmatrix} \begin{pmatrix} \vec{q} \\ \vec{h} \\ \vec{a} \end{pmatrix}. \tag{11}$$

where $\theta$ is the Lorentzian angle between the vectors $\vec{q}^*$ and $\vec{q}$.

**Case 2.** If the Mannheim offset $\varphi^*$ of $\varphi$ is a timelike ruled surface of the type $N_+$ then we have

$$\begin{pmatrix} \vec{q}^* \\ \vec{h}^* \\ \vec{a}^* \end{pmatrix} = \begin{pmatrix} \cosh\theta & \sinh\theta & 0 \\ 0 & 0 & 1 \\ \sinh\theta & \cosh\theta & 0 \end{pmatrix} \begin{pmatrix} \vec{q} \\ \vec{h} \\ \vec{a} \end{pmatrix}. \tag{12}$$

where $\theta$ is the hyperbolic angle between the vectors $\vec{q}^*$ and $\vec{q}$. By definition, the parametrization of $\varphi^*$ can be given by

$$\varphi^*(s,v) = \vec{c}(s) + R(s)\vec{a}(s) + v\vec{q}^*(s). \tag{14}$$

From the definition of $\vec{h}^*$, we get $\vec{h}^* = \dfrac{d\vec{q}^*}{ds} / \left\| \dfrac{d\vec{q}^*}{ds} \right\|$. So that we have $\dfrac{d\vec{q}^*}{ds} = \lambda \vec{h}^*$, ($\lambda$ is a scalar). Using this equality and the fact that the base curve of $\varphi^*$ is also striction curve we get

$\left\langle \dfrac{d}{ds}(\vec{c}+R\vec{a}), \vec{a} \right\rangle = 0$. It follows that $-\left\|\dfrac{d\vec{q}}{ds}\right\| d_\varphi + \dfrac{dR}{ds} = 0$. Thus we can give the following theorems.

**Theorem 4.1.** *Let the timelike ruled surface $\varphi^*$ of the type $N_-$ or $N_+$ be Mannheim offset of the spacelike ruled surface $\varphi$. Then $\varphi$ is a developable spacelike ruled surface if and only if $R$ is a constant.*

Now, we can give the characterizations of the Mannheim offsets of a spacelike ruled surface. Let the timelike ruled surface $\varphi^*$ be Mannheim offset of the developable spacelike ruled surface $\varphi$. By the definition $\varphi^*$ can be of the type $N_+$ or $N_-$. Then we can give the following theorem.

**Theorem 5.1. i)** *Let the timelike ruled surface $\varphi^*$ be the Mannheim offset of the developable spacelike ruled surface $\varphi$. Then $\varphi^*$ is developable if and only if the following equality holds*

*i)* $\cosh\theta - R\kappa \dfrac{ds_1}{ds}\sinh\theta = 0$, *if $\varphi^*$ is of the type $N_-$,* \hfill (15)

*ii)* $\sinh\theta - R\kappa \dfrac{ds_1}{ds}\cosh\theta = 0$, *if $\varphi^*$ is of the type $N_+$.* \hfill (16)

**Proof. i)** Let the timelike ruled surface $\varphi^*$ be of the type $N_-$ and also developable. Then we have

$$\dfrac{d\vec{c}^*}{ds} = \mu \vec{q}^*, \tag{17}$$

where $\mu$ is scalar and $s$ is the arc-length parameter of the striction curve $(c)$ of the spacelike ruled surface $\varphi$. Then from (11) we obtain

$$\dfrac{d\vec{c}}{ds} + \dfrac{dR}{ds}\vec{a} + R\dfrac{ds_1}{ds}\dfrac{d\vec{a}}{ds_1} = \mu(\sinh\theta \vec{q} + \cosh\theta \vec{h}). \tag{18}$$

From Theorem 4.1 and by using (5) we get

$$\vec{q} + R\dfrac{ds_1}{ds}\kappa\vec{h} = \mu\sinh\theta \vec{q} + \mu\cosh\theta \vec{h}. \tag{19}$$

From the last equation it follows that

$$\cosh\theta - R\kappa \dfrac{ds_1}{ds}\sinh\theta = 0. \tag{20}$$

Conversely, if (20) holds then for the tangent vector of the striction curve $(\vec{c}^*)$ of the timelike ruled surface $\varphi^*$ of the type $N_-$ we can write

$$\begin{aligned}
\dfrac{d\vec{c}^*}{ds} &= \dfrac{d}{ds}(\vec{c}+R\vec{a}) \\
&= \vec{q} + R\dfrac{ds_1}{ds}\kappa\vec{h} \\
&= \dfrac{1}{\sinh\theta}(\sinh\theta \vec{q} + \cosh\theta \vec{h}) \\
&= \dfrac{1}{\sinh\theta}\vec{q}^*
\end{aligned}$$

Thus $\varphi^*$ is developable.

**ii)** Let the timelike ruled surface $\varphi^*$ of the type $N_+$ be Mannheim offset of the developable spacelike ruled surface $\varphi$. By making the similar calculations in the proof of the Theorem 4.1 (i) it can be easily shown that $\varphi^*$ is developable if and only if the following equality holds

$$\sinh\theta - R\kappa\frac{ds_1}{ds}\cosh\theta = 0. \tag{21}$$

**Theorem 5.2.** *Let $\varphi$ be a developable spacelike ruled surface. The developable timelike ruled surface $\varphi^*$ of the type $N_+$ or $N_-$ is a Mannheim offset of the spacelike ruled surface $\varphi$ if and only if the following relationship holds*

$$\frac{d\kappa}{ds} = \frac{1}{R}\left(R^2\kappa^2\left(\frac{ds_1}{ds}\right)^2 - 1\right) - \frac{1}{ds_1/ds}\frac{d^2s_1}{ds^2}\kappa. \tag{22}$$

**Proof.** Let the developable timelike ruled surface $\varphi^*$ be a Mannheim offset of the spacelike ruled surface $\varphi$. Assume that $\varphi^*$ is of the type $N_-$. From Theorem 5.1 (i) we have

$$R\kappa\frac{ds_1}{ds} = \coth\theta. \tag{23}$$

Using (11) we have

$$\frac{d\vec{q}^*}{ds} = \cosh\theta\left(\frac{d\theta}{ds} + \frac{ds_1}{ds}\right)\vec{q} + \sinh\theta\left(\frac{d\theta}{ds} + \frac{ds_1}{ds}\right)\vec{h} + \cosh\theta\frac{ds_1}{ds}\kappa\vec{a}. \tag{24}$$

From (24) and definition of $\vec{h}^*$ we have

$$\frac{d\theta}{ds} = -\frac{ds_1}{ds}. \tag{25}$$

Differentiating (23) with respect to $s$ and using (25) we get

$$\frac{d\kappa}{ds} = \frac{1}{R}\left(R^2\kappa^2\left(\frac{ds_1}{ds}\right)^2 - 1\right) - \frac{1}{ds_1/ds}\frac{d^2s_1}{ds^2}\kappa. \tag{26}$$

Conversely, if (26) holds then for nonzero constant scalar $R$ we can define a timelike ruled surface $\varphi^*$ of the type $N_-$ as follows

$$\vec{\varphi}^*(s,v) = \vec{c}^*(s) + v\vec{q}^*(s), \tag{27}$$

where $\vec{c}^*(s) = \vec{c}(s) + R\vec{a}(s)$. Since $\varphi^*$ is developable, we have

$$\frac{d\vec{c}^*}{ds} = \frac{ds^*}{ds}\vec{q}^*, \tag{28}$$

where $s$ and $s^*$ are the arc-length parameters of the striction curves $(\vec{c})$ and $(\vec{c}^*)$, respectively. From (28) we get

$$\frac{ds^*}{ds}\vec{q}^* = \frac{d}{ds}(\vec{c} + R\vec{a}) = \vec{q} + R\kappa\frac{ds_1}{ds}\vec{h}. \tag{29}$$

By taking the derivative of (29) with respect to $s$, we have

$$\frac{d^2s^*}{ds^2}\vec{q}^* + \frac{ds^*}{ds}\frac{d\vec{q}^*}{ds} = R\kappa\left(\frac{ds_1}{ds}\right)^2\vec{q} + \left(\frac{ds_1}{ds} + R\kappa\frac{d^2s_1}{ds^2} + R\frac{ds_1}{ds}\frac{d\kappa}{ds}\right)\vec{h} + R\kappa^2\left(\frac{ds_1}{ds}\right)^2\vec{a}. \tag{30}$$

From the hypothesis and the definition of $\vec{h}^*$, we get

$$\frac{d^2s^*}{ds^2}\vec{q}^* + \frac{ds^*}{ds}\lambda\vec{h}^* = R\kappa\left(\frac{ds_1}{ds}\right)^2\vec{q} + R^2\kappa^2\left(\frac{ds_1}{ds}\right)^3\vec{h} + R\kappa^2\left(\frac{ds_1}{ds}\right)^2\vec{a}, \tag{31}$$

where $\lambda$ is a scalar. By taking the vector product of (29) with (31), we obtain

$$\left(\frac{ds^*}{ds}\right)^2\lambda\vec{a}^* = R^2\kappa^3\left(\frac{ds_1}{ds}\right)^3\vec{q} + R\kappa^2\left(\frac{ds_1}{ds}\right)^2\vec{h}. \tag{32}$$

Taking the vector product of (32) with (29), we have

$$\left(\frac{ds^*}{ds}\right)^3 \lambda \vec{h}^* = \left[R^3\kappa^4\left(\frac{ds_1}{ds}\right)^4 - R\kappa^2\left(\frac{ds_1}{ds}\right)^2\right]\vec{a}. \tag{33}$$

It shows that, the developable timelike ruled surface $\varphi^*$ of the type $N_-$ is a Mannheim offset of the spacelike ruled surface $\varphi$.

If $\varphi^*$ is of the type $N_+$ then making the similar calculations it is easily seen that the developable timelike ruled surface $\varphi^*$ is a Mannheim offset of the spacelike ruled surface $\varphi$ if and only if the same equation i.e. equation (22) holds.

Let now the timelike ruled surface $\varphi^*$ of the type $N_+$ or $N_-$ be a Mannheim offset of the spacelike ruled surface $\varphi$. If the trajectory ruled surfaces generated by the vectors $\vec{h}^*$ and $\vec{a}^*$ of $\varphi^*$ are denoted by $\varphi_{h^*}$ and $\varphi_{a^*}$, respectively, then we can write

$$\vec{q}_1^* = \vec{a}, \; \vec{h}_1^* = \mp\vec{h}, \; \vec{a}_1^* = \mp\vec{q}, \tag{34}$$

$$\vec{q}_2^* = \cosh\theta\vec{q} + \sinh\theta\vec{h}, \; \vec{h}_2^* = \mp\vec{a}, \; \vec{a}_2^* = \pm(\sinh\theta\vec{q} + \cosh\theta\vec{h}), \text{ if } \varphi^* \text{ is of the type } N_-. \tag{35}$$

$$\vec{q}_2^* = \sinh\theta\vec{q} + \cosh\theta\vec{h}, \; \vec{h}_2^* = \mp\vec{a}, \; \vec{a}_2^* = \pm(\cosh\theta\vec{q} + \sinh\theta\vec{h}), \text{ if } \varphi^* \text{ is of the type } N_+. \tag{36}$$

where $\{\vec{q}_1^*, \vec{h}_1^*, \vec{a}_1^*\}$ and $\{\vec{q}_2^*, \vec{h}_2^*, \vec{a}_2^*\}$ are the Frenet Frames of the ruled surfaces $\varphi_{h^*}$ and $\varphi_{a^*}$, respectively. Therefore from (34), (35) and (36) we have the following.

**Corollary 5.3.** (a) $\varphi_{h^*}$ is a Bertrand offset of $\varphi$.

(b) $\varphi_{a^*}$ is a Mannheim offset of $\varphi$.

Now we can give the followings.

Let the timelike ruled surface $\varphi^*$ of the type $N_+$ or $N_-$ be a Mannheim offset of the developable spacelike ruled surface $\varphi$. From (3), (5), (11) and (12), we obtain

$$p_{h^*} = -\frac{1}{(ds_1/ds)\kappa}, \tag{37}$$

$$p_{a^*} = \frac{-\sinh\theta + R(ds_1/ds)\kappa\cosh\theta}{(ds_1/ds)\kappa\sinh\theta}, \text{ if } \varphi^* \text{ is of the type } N_-. \tag{38}$$

$$p_{a^*} = \frac{-\cosh\theta + R(ds_1/ds)\kappa\sinh\theta}{(ds_1/ds)\kappa\cosh\theta}, \text{ if } \varphi^* \text{ is of the type } N_+. \tag{39}$$

Then, we can give the following corollary.

**Corollary 5.4.** (a) $\varphi_{h^*}$ is nondevelopable while $\varphi$ is developable.

(b) $\varphi_{a^*}$ is developable while $\varphi$ is developable if and only if the following equalities holds,

$$-\sinh\theta + R(ds_1/ds)\kappa\cosh\theta = 0, \text{ if } \varphi^* \text{ is of the type } N_-. \tag{40}$$

$$-\cosh\theta + R(ds_1/ds)\kappa\sinh\theta = 0, \text{ if } \varphi^* \text{ is of the type } N_+. \tag{41}$$

**6. Example.** Let consider the spacelike ruled surface $\varphi(s,v)$ given by the parametrization

$$\varphi(s,v) = \left(\cosh(s) + v\frac{\sqrt{2}}{2}\sinh(s), \; v\frac{\sqrt{2}}{2}, \; \sinh(s) + v\frac{\sqrt{2}}{2}\cosh(s)\right),$$

which is rendered in Fig. 1. A Mannheim offset of $\varphi(s,v)$ is

$$\varphi_1^*(s,v) = \left( \cosh(s) - s\frac{\sqrt{2}}{2}\sinh(s) + v\frac{\sqrt{6}}{2}\sinh(s) + 2v\cosh(s),\ s\frac{\sqrt{2}}{2}\sinh^2(s) + v\frac{\sqrt{6}}{2}, \right.$$
$$\left. \sinh(s) - s\frac{\sqrt{2}}{2}\cosh(s) + v\frac{\sqrt{6}}{2}\cosh(s) + 2v\sinh(s) \right)$$

which is a timelike ruled surface of the type $N_-$. (See Fig. 2 (a)).

Another Mannheim offset of $\varphi(s,v)$ is

$$\varphi_2^*(s,v) = \left( \cosh(s) - \frac{3\sqrt{2}}{2}\sinh(s) + v\sqrt{2}\sinh(s) + v\sqrt{3}\cosh(s),\ \frac{3\sqrt{2}}{2}\sinh^2(s) + v\sqrt{2}, \right.$$
$$\left. \sinh(s) - \frac{3\sqrt{2}}{2}\cosh(s) + v\sqrt{2}\cosh(s) + v\sqrt{3}\sinh(s) \right)$$

which is a timelike ruled surface of the type $N_+$. (See Fig. 2 (b)).

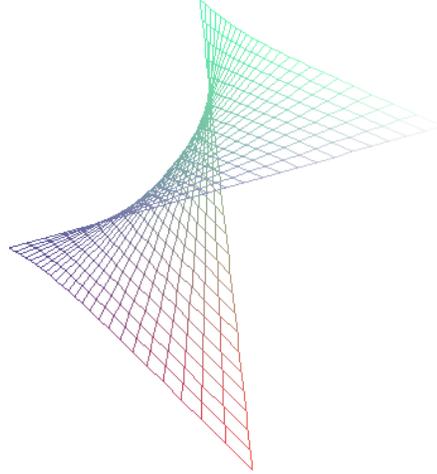

**Fig. 1.** Spacelike ruled surface $\varphi(s,v)$

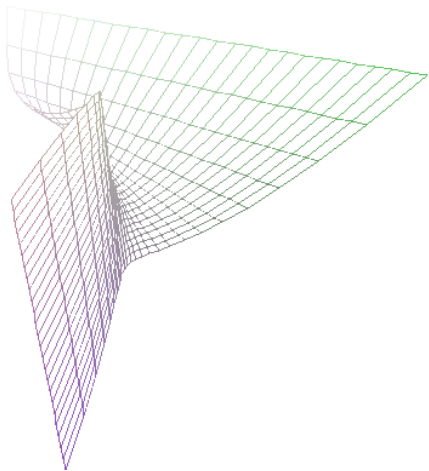 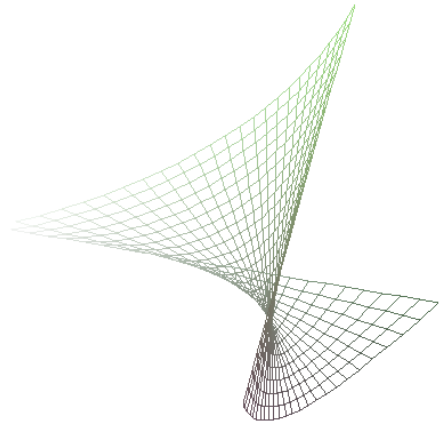

a) Timelike Mannheim offsets $\varphi_1^*(s,v)$ of the type $N_-$

b) Timelike Mannheim offsets $\varphi_2^*(s,v)$ of the type $N_+$

**Fig. 2.**

### 7. Conclusion

In the surface theory, offset surfaces have an important role because of they are used in geometric design problems largely. This paper introduce a new kind of surface offsets called Mannheim offset in Minkowski 3-space. The notion of the Mannheim offsets is given for the spacelike ruled surfaces in Minkowski 3-space $IR_1^3$. It is shown that the Mannheim offsets of a spacelike ruled surface is a timelike ruled surface of the type $N_-$ or $N_+$. Furthermore, it is shown that developable spacelike ruled surfaces can have a developable timelike Mannheim offset if the derivative of the conical curvature $\kappa$ of the directing cone of the spacelike ruled surface holds an equation given by (22).